\documentclass{nhm}
\usepackage{amsmath}
  \usepackage{paralist}
  \usepackage{graphics} 
  \usepackage{epsfig} 
 \usepackage[colorlinks=true]{hyperref}
\hypersetup{urlcolor=blue, citecolor=red}

\def\currentvolume{8}
 \def\currentissue{3}
  \def\currentyear{2013}
   \def\currentmonth{September}
    \def\ppages{627--648}
     \def\DOI{10.3934/nhm.2013.8.627}

\newtheorem{theorem}{Theorem}[section]

\newtheorem{lemma}[theorem]{Lemma}

{\theoremstyle{definition}
\newtheorem{definition}[theorem]{Definition}
\newtheorem{remark}{Remark}
}

\def\sqr#1#2{\vbox{\hrule height .#2pt
\hbox{\vrule width .#2pt height #1pt \kern #1pt
\vrule width .#2pt}\hrule height .#2pt }}
\def\square{\sqr74}
\def\endproof{\hphantom{MM}\hfill\llap{$\square$}\goodbreak}

\def\begi{\begin{itemize}}
\def\endi{\end{itemize}}
\def\bega{\begin{array}}
\def\enda{\end{array}}

\def\forall{\hbox{for every~ }}
\def\ov{\overline}
\def\ds{\displaystyle}
\def\R{\mathbb{R}}

\def\L{{\bf L}}

\def\B{{\mathcal B}}

\def\forall{\hbox{for all}~}

\def\vp{\varphi}

\def\I{{\mathcal I}}
\def\O{{\mathcal O}}
\def\G{{\mathcal G}}

\def\U{{\mathcal U}}

\def\C{{\mathcal C}}

\def\c{\centerline}

\def\wto{\rightharpoonup}

\def\v{\vskip 1em}

\def\ve{\varepsilon}

\def\V{{\mathcal V}}

\def\bel{\begin{equation}\label}
\def\eeq{\end{equation}}

\title[Optima and Equilibria on Traffic Networks]
      {Existence of Optima and Equilibria for Traffic Flow on Networks}

\author[Alberto Bressan and Ke Han]{}


 \email{bressan@math.psu.edu}
 \email{kxh323@psu.edu}

\thanks{This work is partially supported by NSF through grant EFRI-1024707, ÔÔA theory of complex transportation network designÕÕ}

\begin{document}
\maketitle

\centerline{\scshape Alberto Bressan and Ke Han}
\medskip
{\footnotesize
 \centerline{Department of Mathematics}
   \centerline{The Pennsylvania State University}
   \centerline{University Park, PA 16802, USA}
} 

%
%
%

\begin{abstract}
This paper is concerned with a
conservation law model of traffic flow on a network of roads,
where each driver chooses his own departure time in order
to minimize the sum of a departure cost and an arrival cost.
The model includes
various groups of drivers,  with different origins and destinations
and having different cost functions.  Under a natural set of assumptions,
two main results are proved: (i)
the existence of a globally optimal solution, minimizing the sum of the costs
to all drivers, and (ii) the existence of
a Nash equilibrium solution, where no driver can lower his own cost by changing
his departure time or the route taken to reach destination.
In the case of Nash solutions, all departure rates are uniformly bounded and have
compact support.
\end{abstract}

\section{Introduction}

System-optimal and user-optimal traffic flows are the two primary objects of study in {\it dynamic traffic assignment} (DTA); they correspond to the dynamic extension of Wardrop's second and first principles, respectively \cite{Wardrop}. These two types of traffic assignment problems are essential for understanding congestion externalities, network performance and the price of anarchy. The  {\it system-optimal} dynamic traffic assignment (SO-DTA) problem aims at minimizing the total travel cost on a system level subject to the constraints of travel demand, link dynamics, flow propagation and travel delay. The user-optimal problem, also known as the {\it dynamic user equilibrium} (DUE), requires that unit travel cost is identical for those route and departure time choices selected by travelers between a given origin-destination pair, and that no road user can
reduce his own cost under such equilibrium flow.

This paper is concerned with the existence of continuous-time system-optimal
and user-optimal traffic flows on a road network, with a specific type of traffic
flow model. In the authors' earlier work \cite{BH1, BH2}, the existence and
uniqueness
of system-optimal and user-optimal solutions were established for a single link,
when the Lighthill-Whitham-Richards model \cite{LW, R} is employed.
In this paper, we extend the existence results to a network of roads.
In particular, the existence
of a system-optimal solution is obtained by taking the limit of a minimizing
sequence of admissible departure distributions. The existence result for user
equilibrium solutions is proved by a fixed point argument, after establishing
a uniform upper bound on all departure rates.         This topological
technique allows  great  generality but
does not yield information about uniqueness or continuous
dependence of the user equilibrium.

Consider a model of traffic flow where  drivers
travel on a network of roads.
We denote by $A_1,\ldots, A_m$ the nodes of the network,
and by $\gamma_{ij}$ the arc connecting $A_i$ with $A_j$.
Following the classical papers \cite{LW, R},
along  each arc the flow of  traffic will be modeled by the conservation law
\bel{clawij}\rho_t + [\rho\, v_{ij}(\rho)]_x~ =~ 0\,.\eeq
Here $t$ is time and $x\in [0, L_{ij}]$ is the space variable along the arc
$\gamma_{ij}$.
By $\rho=\rho(t,x)$ we denote  the traffic density, while the map $\rho\mapsto v_{ij}(\rho)$ is
the speed of cars as function of the density, along the arc $\gamma_{ij}$.
We assume that $v_{ij}$ is a continuous, nonincreasing function.
If $v_{ij}(0) >0$ we say that the arc $\gamma_{ij}$ is {\it viable}.
It is quite possible that
two nodes $i,j$ are not directly linked by a road.   This situation can be easily modeled by taking
$v_{ij}\equiv 0$, so that the arc is not viable.
The conservation laws (\ref{clawij}) are supplemented by
suitable boundary conditions at
points of junctions, which will be discussed later.

We consider $n$ groups of drivers traveling on the network.
Different groups are distinguished by the locations of departure and arrival, and
by their cost functions.
For $k\in \{1,\ldots,n\}$, let  $G_k$ be the total number of drivers in the $k$-th group.
All these drivers
depart from a node $A_{d(k)}$ and arrive at a node $A_{a(k)}$, but can choose
different paths to reach destination.
Of course, we assume that there exists at least one chain of viable arcs
\bel{path}
\Gamma~\doteq~ \Big(\gamma_{\strut i(0), i(1)}\,, ~\gamma_{\strut i(1), i(2)}\,,~\ldots~,\,
 \gamma_{\strut i(\nu-1), i(\nu)}
\Big)\eeq
with $i(0)= d(k)$ and $i(\nu)= a(k)$,
connecting the departure node $A_{d(k)}$ with the arrival node $A_{a(k)}$.
We shall denote by
$$\V~\doteq~\Big\{ \Gamma_1,\, \Gamma_2,~\ldots~,\, \Gamma_N\Big\}$$
the set of all viable chains (i.e.~concatenations of viable arcs)
which do not contain any  closed loop.
Since there are $m$ nodes, and each chain can visit each of them at most once,
the cardinality of $\V$ is bounded by $(m+1)!$.
For a given $k\in \{1,\ldots,n\}$,
let $\V_k\subset\V$ be the set of all viable paths for the $k$-drivers,
connecting $A_{d(k)}$ with $A_{a(k)}$.

By $U_{k,p}(t)$ we denote the total number of drivers of the $k$-th group, traveling
along the viable path $\Gamma_p$,
who  have started their journey before time $t$.
 \v
\begin{definition}\label{definition1}
A {\bf departure distribution
function} $t\mapsto U_{k,p}(t)$
is  a bounded, nondecreasing, left-continuous function,
such that
$$U_{k,p}(-\infty)~\doteq~ \lim_{t\to -\infty} U_{k,p}(t) ~=~0\,.$$
Given group sizes $G_1,\ldots, G_n\geq 0$, we say that
a set of departure distribution functions $\{ U_{k,p}\}$
is {\bf admissible}
if it satisfies the constraints
\bel{c1} \sum_p  U_{k,p}(+\infty) ~=~G_k\qquad\qquad \qquad k=1,\ldots, n\,.\eeq
\end{definition}

\v
Since $G_k$ is the total number of drivers in the $k$-th group,
the admissibility condition (\ref{c1}) means that, sooner or later,
every driver of each group has to depart.
If the function $U_{k,p}$ is absolutely continuous, its derivative will be denoted by
\bel{Udef1}
\ov u_{k,p}(t)~ =~ {d\over dt} U_{k,p}(t)\,.\eeq
Clearly, $\bar u_{k,p}$ measures the rate of  departures
 of $k$-drivers traveling along $\Gamma_p$.

The overall traffic pattern can be determined by (i) the departure distribution functions
$U_{k,p}(\cdot)$, (ii) the conservation laws
(\ref{clawij}), and (iii) a suitable set of conditions at junctions, specifying the priorities assigned to drivers that wish to enter the same road.

In this paper we consider the simplest type of condition
at junctions,  where a separate queue can
form at the entrance of each road.
Drivers arriving at the node $A_i$
from all incoming roads,
and who want to travel along the arc $\gamma_{ij}$,
join a queue at the entrance of this outgoing arc.  Their place in the queue
is determined by the time at which they arrive at $A_i$, first-come first-serve.
Some additional care is needed to handle the case where different groups of drivers
depart from the same node. Indeed, if a positive amount of drivers initiate their
journey exactly at the same time, some additional information is needed to determine
their place in the queue.   This can be
achieved in terms of the {\it prioritizing functions} introduced in \cite{BH2}.

As in \cite{BH1, BH2, BLSY},
we consider a set of departure costs $\vp_k(\cdot)$, and arrival costs $\psi_k(\cdot)$
for the various drivers.
Namely, a driver of the $k$-th group departing at time $\tau^d$ and arriving
at destination at time $\tau^a$ will incur in the total cost
\bel{ck}
\vp_k(\tau^d) +\psi_k(\tau^a).\eeq
In this framework, the concepts of globally optimal solution and of
Nash equilibrium solution considered in \cite{BH1, BH2}
can  be extended to traffic flows on a network of roads.
\v
\begin{definition}\label{definition2}    An admissible family
$\{ U_{k,p}\}$ of departure distributions is {\bf globally optimal}
if it minimizes the sum of the total costs of all drivers.
\end{definition}

\v
\begin{definition}\label{definition3}
An admissible family
$\{ U_{k,p}\}$ of departure distributions is a
{\bf Nash equilibrium solution}
if no driver of any group can lower his own total  cost by
changing departure time or switching to a different path to reach destination.
\end{definition}
\v

The main goal of
this paper is to prove the existence of a globally optimal solution and of a Nash
equilibrium solution, under natural assumptions on the costs and on the flux functions.

In the case of a single group of drivers traveling on a single road,
the existence and uniqueness of such solutions were proved in
\cite{BH1}.
We highlight the main features of the present analysis.
\begi
\item Following the direct method of the Calculus of Variations, a
globally optimal solution is constructed by taking the limit of a minimizing
sequence $\{ U^{(\nu)}_{k,p}\}_{\nu\geq 1}$ of admissible departure distributions.
The existence of the limit is
guaranteed by the ``tightness" of the sequence of approximating measures.
Namely, for each $\ve>0$ there exists $T>0$ (independent of $n$) such that the
total amount of drivers departing at times $t\notin [-T,T]$ is less than $\ve$.

\item Toward the existence of a Nash equilibrium, the
results in \cite{BH1, BH2} used the assumption $\psi_k'\geq 0$, meaning that the arrival cost functions are nondecreasing.
We now strengthen this assumption to $\psi_k'>0$, so that the
arrival costs are strictly increasing.  This apparently minor change in the hypotheses
has an important consequence.  Namely, it allows us to prove
a crucial a priori bound on all departure rates,
in any Nash equilibrium solution.

\item The proofs in  \cite{BH1, BH2} relied on a monotonicity argument.   Indeed, the
departure distribution $U(\cdot)$ for a Nash equilibrium was obtained as the (unique)
pointwise supremum of
a family of admissible distributions, satisfying an additional constraint.
On the other hand, the present existence result is proved by a
fixed point argument.   By its nature, this topological  technique cannot  yield information about uniqueness or continuous  dependence of the Nash equilibrium.
\endi

The paper is organized as follows.  In Section 2 we describe more carefully
the traffic flow model, explaining how to compute the admissible solutions
using the Lax-Hopf formula.
In Section~3 we prove the existence of a globally optimal solution, while
Section~4 is devoted to the existence of a Nash equilibrium solution.

For the modeling of traffic flow we refer to \cite{Bellomo, Dag, LW, R}.
Traffic flow on networks has been the topic of an extensive
literature, see for example \cite{CGP,  F, FKKR, GP} and references therein.
A different type of optimization problems for traffic flow
was considered in \cite{GHKL}.
\v
\section{Analysis of the traffic flow model}

In our  model, $x\in [0,L_{ij}]$ is the space variable, describing
a point along the arc $\gamma_{ij}$.  Here $L_{ij}$
measures the length of this arc.
The basic assumptions on the flux functions $F_{ij}(\rho) = \rho\, v_{ij}(\rho)$ and
on the cost functions $\vp_k,\psi_k$
are as follows.

\begi
\item[{\bf (A1)}] For every viable arc $\gamma_{ij}$, the flux function $\rho\mapsto
F_{ij}(\rho) = \rho v_{ij}(\rho)$ is  continuous, concave down, and non-negative
on some interval $[0, \ov  \rho_{ij}]$, with  $F_{ij}(0) = F_{ij}(\ov \rho_{ij}) =0$.
We shall denote by $\rho_{ij}^*\in ]0, \ov\rho_{ij}[\,$ the unique value such that
\bel{fprop} F(\rho_{ij}^*) ~=~F^{max}_{ij}
\doteq~\max_{\rho\in [0, \ov \rho_{ij}]} F_{ij}(\rho)\,,\qquad\qquad
F'_{ij}(\rho) ~>~0\qquad \hbox{for a.e.}~~\rho\in [0, \rho_{ij}^*]\,.
\eeq

\item[{\bf (A2)}] For every $k\in \{1,\ldots,n\}$
the cost functions $\vp_k,\psi_k$ are continuously differentiable
and satisfy
\bel{costp}\left\{
\bega{l}\vp'_k(t)~<~0\,,\cr \cr
 \psi'_k(t)~> ~0\,,\enda\right.\qquad\qquad
\lim_{|t|\to\infty}
\Big(\vp_k(t)+ \psi_k(t)\Big)~=~+\infty\,. \eeq
\endi
\v

\begin{figure}[htbp]
   \centering
 \includegraphics[width=1.0\textwidth]{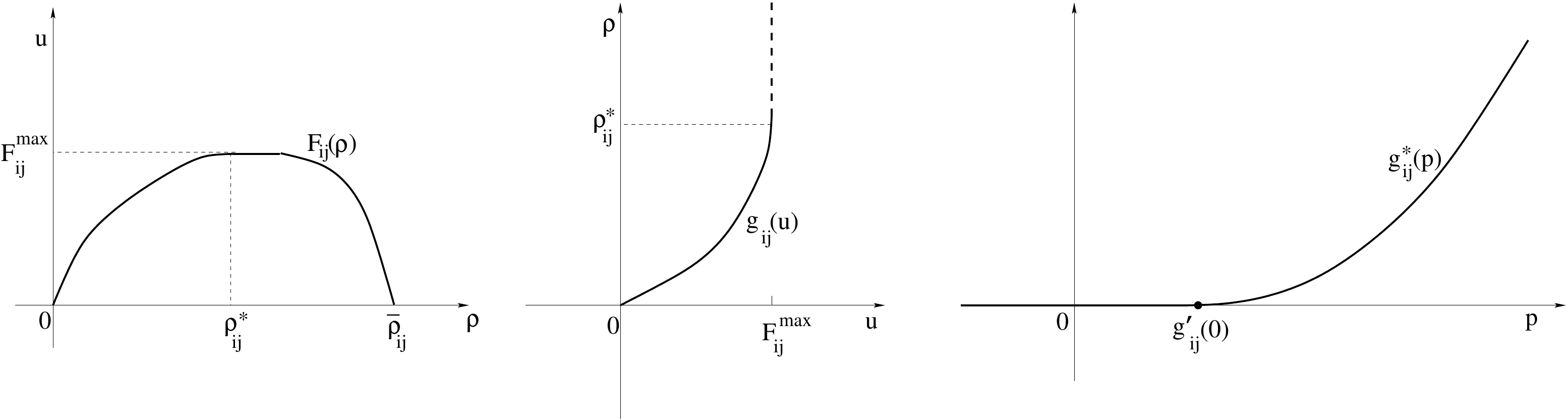}
   \caption{\small  Left: the function $\rho\mapsto F_{ij}(\rho)=
   \rho\,v_{ij}(\rho)$ describing the flux of cars. Middle:
the  function $g_{ij}$, implicitly defined
by $g_{ij}(\rho v_{ij}(\rho)) = \rho$.
Right: the Legendre transform
   $g_{ij}^*$.}
   \label{f:t94}
\end{figure}

\begin{remark}\label{rm1}
 {\rm By (A1), the flux function $u=F_{ij}(\rho)$ is continuous, concave
and strictly
increasing on the interval $[0, \rho^*_{ij}]$.   Therefore, it has
a continuous inverse: $\rho = g_{ij}(u)$. As shown in fig.~\ref{f:t94},
the function $u\mapsto g_{ij}(u)$ is
convex and
maps the interval $[0, F_{ij}^{max}]$ onto $[0, \rho^*_{ij}]$.
The concavity of $F_{ij}$ means that
$$F_{ij}(\theta \rho + (1-\theta)\rho')~\geq~\theta F_{ij}( \rho)
+ (1-\theta)F_{ij}(\rho')\qquad\qquad \theta\in
[0,1],~~~\rho,\rho'\in [0, \bar\rho_{ij}]\,.$$
}
\end{remark}

\begin{remark} {\rm
According to (A2), the cost for early departure is strictly
decreasing in time, while cost for late arrival is
strictly increasing. The assumption that these costs tend to infinity as
$t\rightarrow\pm \infty$ coincides with common sense and guarantees that
in an equilibrium  solution the departure rates are compactly supported.}
\end{remark}

\begin{remark} {\rm
In the engineering literature (see for example \cite{FKKR})
it is common to define the travel
cost as the  sum of the travel time plus a penalty if the arrival time does not coincide with the target
time $T_A$ :
\bel{engrcost}
D(t)+\Psi\big(t+D(t)-T_A\big).
\eeq
Here $D(t)$ denotes the total duration of the trip for a driver departing at
time $t$,
while $\Psi$ is a penalty function.
Calling $\tau^a(t)= t+D(t)$ the arrival time of driver departing at $t$, the cost function (\ref{engrcost})
can be recast in the form (\ref{ck}).  Indeed,
\begin{align*}
D(t)+\Psi\big(t+D(t)-T_A\big)&~=~-t+t +D(t)+\Psi\big(t+D(t)-T_A\big)\\
&~=~-t+\tau^a(t)+\Psi\big(\tau^a(t)-T_A\big)\\
&~=~\vp(t)+\psi(\tau^a(t))\,,
\end{align*}
where $\vp(t)\doteq -t$, $\psi(\tau)\doteq \tau+\Psi\big(\tau-T_A)$.
In order that the assumptions {\bf (A2)} be satisfied, it suffices to
require that the function $\Psi$  be
continuously differentiable and
\bel{engramp}
 \Psi~\geq~0,\qquad\qquad
\Psi'~>~-1,\qquad \qquad\lim_{|t|\rightarrow \infty} (\Psi(t)-t)~=~+\infty\,.
\eeq
}
\end{remark}

\v
\subsection{Traffic flow with an absolutely continuous departure distribution.}
We now describe more in detail how the traffic flow on the entire network
can be uniquely determined, given the departure distributions $U_{k,p}$.
As a first step, we consider the absolutely continuous case,
so that (\ref{Udef1}) holds.

Along any arc $\gamma_{ij}$, the traffic density $\rho_{ij}$ satisfies
a boundary value problem of the form
\bel{bvp}
\begin{cases}
\partial_t \rho_{ij}(t,\,x)_t+\partial_xF_{ij}\big(\rho_{ij}(t,\,x)\big)~=~0
\qquad &(t,\,x)\in \R\times [0,\,L_{ij}]\,, \\
F_{ij}\big(\rho_{ij}(t,\,0)\big)~=~ u^-_{ij}(t)\qquad & t\in \R\,,
\end{cases}
\eeq
where $u_{ij}^-(t)$ describes the incoming flux at $x=0$.
Following the approach in \cite{BH1, BH2}, we switch the roles of the variables
$t,\,x$, replacing the above boundary value problem (\ref{bvp}) with a Cauchy
problem for the conservation law describing the flux $u_{ij}=F_{ij}(\rho_{ij})$:
\bel{ivp}
\begin{cases}
\partial_x u_{ij}(t,\,x) + \partial_t g_{ij}\big(u_{ij}(t,\,x)\big)~=~0\qquad &(t,\,x)\in \R\times [0,\,L_{ij}]\,,  \\
u_{ij}(t,\,0)~=~u_{ij}^-(t)\qquad & t\in\R\,.
\end{cases}
\eeq
Here $g_{ij}:[0, F_{ij}^{max}] \mapsto [0, \rho_{ij}^*]$ is
the inverse of the function $F_{ij}$, as in Remark~\ref{rm1}.
Consider the integrated functions
$$
U_{ij}(t,\,x)~=~\int_{-\infty}^t u_{ij}(s,\,x)\,ds,
\qquad U_{ij}^-(t)~=~\int_{-\infty}^t  u_{ij}^-(s)\,ds\,.
$$
Then $U_{ij}(t,\,x)$ provides a solution to the  Hamilton-Jacobi equation
\bel{ivpHJ}
\begin{cases}
\partial_x U_{ij}(t,\,x)+g_{ij}\big(\partial_t U_{ij}(t,\,x)\big)~=~0\,,\\
U_{ij}(t,\,0)~=~ U_{ij}^-(t)\,.
\end{cases}
\eeq
The viscosity solution to the above Cauchy problem is given by the Lax-Hopf formula \cite{Evans}
\bel{Lax}
U_{ij}(t,\,x)~=~\min_{\tau}\left\{ U_{ij}^-(\tau)+x\,g_{ij}^*\Big({t-\tau\over x}\Big)\right\},\qquad x\in[0,\,L_{ij}]\,.
\eeq
Here $g_{ij}^*$ denotes the Legendre transform of $g_{ij}$, namely
\bel{g*def}
g_{ij}^*(p)~\doteq~\max_{u\in [0,\, F_{ij}^{max}]}~\Big( pu - g_{ij}(u)\Big).
\eeq

\begin{figure}[htbp]
   \centering
 \includegraphics[width=0.5\textwidth]{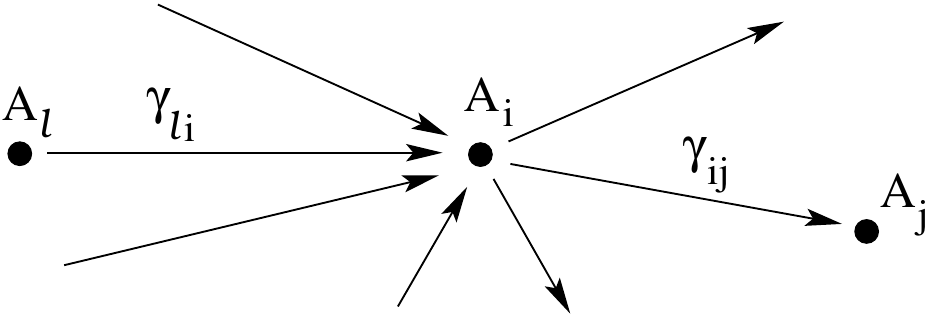}
   \caption{\small  A generic node of the network.}
   \label{f:t95}
\end{figure}

For any given node $A_i$, in general
there will be several incoming arcs $\gamma_{\ell,i}$, $\ell\in
\I(i)$ and several outgoing arcs $\gamma_{i,j}$, $j \in \O(i)$  (see Figure  \ref{f:t95}).
To determine the flux at the entrance of the arc $\gamma_{ij}$,
one needs
to know how many drivers, after reaching the node $A_i$,
actually want to take the road $\gamma_{ij}$.  For this purpose, we need to introduce
the distribution functions
\begi
\item{}
$U_{kp, ij}^-(t)~\doteq$~ {\it total number of drivers of the $k$-th group, traveling along the path $\Gamma_p$, that have entered the arc $\gamma_{ij}$ (possibly joining a queue at the entrance) within time $t$.}
\item{}
$U_{kp, ij}^+(t)~\doteq$~ {\it total number of drivers of the $k$-th group, traveling along the path $\Gamma_p$, that have exited from the arc $\gamma_{ij}$ (reaching the
node $A_j$) within time $t$.}
\endi
Of course, $U_{kp,ij}^\pm\equiv 0$ if the path $\Gamma_p$ does not contain the arc $\gamma_{ij}$.

The entire flux along the network is entirely determined by the
functions $U_{kp,ij}^\pm$.   To recursively compute these functions,
observe that the total number of drivers who have entered the arc $\gamma_{ij}$
(possibly joining a queue) within time $t$ is determined by
\bel{Uij}U_{ij}^-(t)~=~\sum_{\ell\in \I(i)} \sum_{k}\sum_{\gamma_{\ell i},
\gamma_{ij}\in \Gamma_p}
U_{kp, \ell i}^+(t) ~+~\sum_{d(k)=i}\sum_{ \gamma_{ij}\in \Gamma_p}~ U_{kp}(t)\,.\eeq
Notice that the first sum accounts for the drivers that transit through the node
$A_i$, while the second sum accounts for the drivers that initiate their journey from the node $A_i$.
Using the solution formula (\ref{Lax}) with $x=L_{ij}$ we obtain
\bel{Uij+}
U_{ij}^+(t)~\doteq U_{ij}(t, L_{ij})~=
~\min_{\tau}\left\{ U_{ij}^-(\tau)+L_{ij}\,g_{ij}^*
\Big({t-\tau\over L_{ij}}\Big)\right\}\,.
\eeq
Notice that the
function $U_{ij}^+$ is nondecreasing and Lipschitz continuous.
Indeed, the rate at which
drivers arrive at the end of the arc $\gamma_{ij}$ is computed by
$$u^+_{ij}(t) ~=~{d\over dt} U^+_{ij}(t)~\in~[0, \, F_{ij}^{max}]\qquad\qquad
\hbox{for a.e.}~~t\,.$$
Among the drivers who reach the end of the arc $\gamma_{ij}$ within time $t$,
we still need  to compute how many belong to the various groups.

In the case where all departure rates $U_{k,p}$
with $d(k) =i$ are absolutely continuous, the functions $U^+_{kp,ij}$ can be
computed as follows.   For a.e.~$t$   there exists a unique time $   \tau^{enter}(t)$
such that
$$U^-_{ij}(   \tau^{enter}(t)) ~=~U^+_{ij}(t).$$
A  driver entering the arc $\gamma_{ij}$ at time $   \tau^{enter}(t)$
thus reaches the end of the arc at time $t$.
The first-come first-serve assumption now implies
\bel{ukp1}
u_{kp, ij}^+(t) ~=~{ u_{kp, \ell i}^+(   \tau^{enter}(t))\over u_{ij}^-(   \tau^{enter}(t))}
\qquad\qquad \hbox{if}\quad \ell\in \I(i)\,,~~\gamma_{\ell i},\gamma_{ij}
\in\Gamma_p\,,\eeq
\bel{ukp2}
u_{kp, ij}^+(t) ~=~{ u_{kp}(   \tau^{enter}(t))\over u_{ij}^-( \tau^{enter}(t))}
\qquad\qquad \hbox{if}\quad \gamma_{ij}
\in\Gamma_p\,,~~i=d(k)\,.\eeq
In other words, among all drivers
 who exit from the arc $\gamma_{ij}$ at time $t$,
 the percentage of $(k,p)$-drivers
must be equal to the percentage of $(k,p)$-drivers that
enter $\gamma_{ij}$ at time $\tau^{enter}(t)$.
In turn, the arrival distributions are computed by
\bel{U99}U^+_{kp, ij}(t)~=~\int_{-\infty}^t u^+_{kp,ij}(t)\, dt\,.\eeq

\begin{remark} {\rm  The above equations  can be
solved iteratively in time. Namely, let
$$\Delta_{min}~\doteq~\min_{ij} {L_{ij}\over v_{ij}(0)}$$
be the minimum time needed to travel along any viable arc of the network.
Given the departure rates $\bar u_{k,p}$,
if the functions $U_{kp,ij}^\pm(t)$ are known for all $t\leq \tau$,
by the above equations (\ref{Uij})--(\ref{U99}) one can uniquely determine the
values of $U_{kp,ij}^\pm(t)$ also for $t\leq \tau + \Delta_{min}$.
}\end{remark}

Let $G_{k,p}$ be the total number of drivers of the $k$-th group who travel
along the path $\Gamma_p$.   The admissibility condition implies
$G_{k,1}+\cdots+ G_{k, N} = G_k$.
We use the Lagrangian variable
$\beta\in [0,G_{k,p}]$ to label a particular driver in the subgroup $\G_{k,p}$
of $k$-drivers traveling along the path $\Gamma_p$.
The departure and arrival time of this
driver will be denoted by
$\tau^d_{k,p}(\beta)$ and $\tau^a_{k,p}(\beta)$, respectively.
Let $U_{k,p}^{depart}(t)=U_{k,p}(t)$  denote the
amount of drivers of the subgroup $\G_{k,p}$ who have departed
before time $t$.   Similarly, let $U^{arrive}_{k,p}(t)$
be the amount of $(k,p)$-drivers who have arrived at destination before time $t$.
For a.e.~$\beta\in [0, G_{k,p}]$
we then have
\bel{tda}
\tau_{k,p}^d(\beta)~=~\inf \Big\{ \tau\,;~~U_{k,p}^{depart}(t)~\geq~\beta\Big\}\,,
\qquad\qquad
\tau_{k,p}^a(\beta)~=~\inf \Big\{ \tau\,;~~U_{k,p}^{arrive}(t)~\geq~\beta\Big\}\,.
\eeq
With this notation, the definition of globally optimal and of
Nash equilibrium solution can be
more precisely formulated.
\v
\begin{definition}\label{definition2'}  An admissible family of departure distributions
$\{U_{k,p}\}$
is a {\bf globally optimal solution} if it provides a global minimum to the
functional
\bel{JJdef}
J~\doteq~\sum_{k,p} \int_0^{G_{k,p}} \Big(\vp_k(\tau^d_{k,p}(\beta)) + \psi_k
(\tau^a_{k,p}(\beta))\Big) d\beta\,.\eeq
\end{definition}

\begin{definition}\label{definition3'} An admissible family of departure distributions
$\{U_{k,p}\}$
is a {\bf Nash equilibrium solution} if
there exist  constants $c_1,\ldots, c_n$ such that:
\begin{itemize}
\item[(i)] For almost every $\beta\in [0,G_{k,p}]$ one has
\bel{same}\vp_k(\tau^q_{k,p}(\beta)) + \psi_k(\tau^a_{k,p}(\beta)) ~=~c_k\,.\eeq
\item[(ii)] For all $\tau\in\R$, there holds
\bel{nobetter}
\vp_k(\tau) + \psi_k(A_{k,p}(\tau))~\geq~c_k\,.\eeq
\end{itemize}
Here $A_{k,p}(\tau)$ is the arrival time of a driver that starts at time $\tau$
from the node $A_{d(k)}$ and
reaches the node $A_{a(k)}$ traveling along the path $\Gamma_p$.
\end{definition}
\v
In other words, condition (i) states that
all $k$-drivers bear the same cost~$c_k$.
Condition (ii) means that, regardless of the
starting time $x$, no $k$-driver can achieve a cost
$<c_k$.
\v
\subsection{Traffic flow with general departure distribution}

If all departure distributions $U_{k,p}$ are absolutely continuous,
the previous analysis shows that
the first-come first-serve assumption on the queues
completely determines the traffic pattern.  This is no longer true if
a positive amount of drivers of different groups initiate their journey
exactly at the same time.   For example, assume that
$i = d(k) = d(k')$ is the departure node for both $k$-drivers and
$k'$-drivers. Let $\gamma_{ij}$ be the first arc in the paths $\Gamma_p$ and
$\Gamma_{p'}$, and assume that,
at the instant $t_0$,
a positive amount of drivers in the subgroups $\G_{k,p}$ and $\G_{k', p'}$
initiate their journey.   Since all these drivers join the
queue at the entrance of the arc $\gamma_{ij}$ at the same time,
additional information must be provided to determine their
relative position in the queue.
For this purpose, we follow the approach introduced in \cite{BH2}.

Given an arc $\gamma_{ij}$, call
\bel{Gij}
\G_{ij}~=~\bigcup_{d(k)=i,\, \gamma_{ij}\in \Gamma_p} \G_{k,p}\eeq
the family of all drivers that initiate their journey from
the node $A_i$, traveling along $\gamma_{ij}$ as first leg of their journey.
The total number of these drivers is
$$G_{ij} ~=~\sum_{d(k)=i,\, \gamma_{ij}\in \Gamma_p} G_{k,p}\,.$$
The {\it cumulative departure distribution}
$U^{depart}_{ij}:\R\mapsto [0, G_{ij}]$
 is defined as
\bel{cumstart}
U^{depart}_{ij}(t)~\doteq~\sum_{d(k)=i,\, \gamma_{ij}\in \Gamma_p} U_{k,p}(t)\,.\eeq
In other words, $U^{depart}_{ij}(t)$ is the total number of drivers in the family
$\G_{ij}$
that depart before time $t$.
The relative position of these drivers in the
queue at the entrance of the arc $\gamma_{ij}$ will be
determined by a set of prioritizing functions.
\v
\begin{definition}\label{definition4}   A set of  {\bf prioritizing functions} for departures
on the arc $\gamma_{ij}$ is a family of nondecreasing maps
\bel{priorf}\B_{k,p}:[0,G_{ij}] ~\mapsto~[0,G_{k,p}]\,,\qquad\qquad
 \G_{k,p}\subseteq\G_{ij}\,,\eeq
satisfying
\bel{pf1}
\sum_{\G_{k,p}\subseteq\G_{ij}}~ \B_{k,p}(\beta) ~=~\beta\qquad\qquad
\forall \beta\in [0,G_{ij}]\,,\eeq
\bel{pf2}
\B_{k,p}(U^{depart}_{ij}(t))~=~U_{k,p}(t)\qquad\qquad\hbox{for a.e.}~~t\in\R\,.
\eeq
\end{definition}
\v
Otherwise stated,
among the first $\beta$ drivers that depart and choose $\gamma_{ij}$ as the first
arc of their journey,
 $\B_{k,p}(\beta)$ counts how many belong to the subgroup $\G_{k,p}$.
If all functions $U_{k,p}$ are continuous, then the
conditions (\ref{pf1})-(\ref{pf2}) uniquely determine the prioritizing
functions $\B_{k,p}$.
In this case, there is actually no need to introduce this concept.
On the other hand, if two or more functions $U_{k,p}$ have an upward jump
at a time $t_0$, then different prioritizing functions are possible.
In practice, this models a situation where a large number of drivers,
of different groups, start their journey by getting to the
entrance of a highway all at the same time.
Their relative positions in an entrance
queue can then be regulated by a traffic light signal, which assigns different priorities to drivers of different groups.

We now show that, by choosing one set of prioritizing functions, the traffic flow on the entire network is entirely determined.
Indeed, for all $k,p,i,j$ the values $U_{kp, ij}^+(t)$, determining how many drivers
of the subgroup $\G_{k,p}$ reach the end of the arc $\gamma_{ij}$ within time $t$,
are computed as follows.
\v
\noindent\textbf{Case  1.}  The only drivers traveling on the arc $\gamma_{ij}$ are those who started
their journey from the node $A_i$.
In this case, as in (\ref{Uij+}) the total number of drivers
arriving at the end of the arc $\gamma_{ij}$
within time $t$ is
\bel{Uij++}U_{ij}^+ (t)~=~
\sum_{k,p} U_{kp, ij}^+(t)~=~\min_{\tau}\left\{ U_{ij}^{depart}(\tau)+L_{ij}\,g_{ij}^*
\Big({t-\tau\over L_{ij}}\Big)\right\}\,.\eeq
Given the prioritizing functions $\B_{k,p}$, the values
$U_{kp, ij}^+(t)$ are immediately obtained by the formula
\bel{U4}
U_{kp, ij}^+(t)~=~\B_{k,p}\left(U_{ij}^+(t)\right)\,.\eeq
\v
\noindent\textbf{Case  2.} The arc $\gamma_{ij}$ is also traveled by drivers who transit through the node
$A_i$, departing from other nodes.  In this case,
drivers originating from $A_i$ have to merge with drivers in transit,
coming from other nodes.
The cumulative distribution function, accounting for the total number
of drivers that have entered the arc $\gamma_{ij}$ within time $t$ is
\bel{incj}U_{ij}^-(t)~=~U_{ij}^{depart}(t) + U^{transit}_{ij}(t)~
=~\sum_{d(k)=i,\gamma_{ij}\in \Gamma_p} U^-_{kp, ij}(t) + \sum_{d(k)\not=i,\gamma_{ij}\in \Gamma_p} U^-_{kp, ij}(t)\,.\eeq
As before, the total number of drivers who have exited from the arc $\gamma_{ij}$
before time $t$ is given by (\ref{Uij+}). To determine how many of these drivers belong to each subgroup
$\G_{k,p}$, we proceed as follows.

Consider the driver who exits from the arc $\gamma_{ij}$ at time $t$.
This driver will have entered the arc $\gamma_{ij}$
at an earlier time $\tau=\tau^{enter}(t) $, such that
\bel{tdep}
\lim_{s\to \tau-}U_{ij}^-(s)~\leq~U^+_{ij}(t)~\leq ~\lim_{s\to \tau+}U_{ij}^-(s)\,.\eeq
Notice that $\tau^{enter}(t)$ is uniquely determined, for all but countably
many times $t$.    For all subgroups of drivers in transit, the first-come first-serve
priority assumption implies
\bel{pr4}
U^+_{kp, ij}(t) ~=~U^-_{kp,ij}(\tau^{enter}(t)).\eeq
Here the right hand side of (\ref{pr4}) is uniquely determined,
because the distribution function
$U^-_{kp, ij}$ in (\ref{incj}) is Lipschitz continuous.

On the other hand, for drivers who initiate their journey at $A_i$, the
distribution function
$U^-_{kp, ij}$ can have a jump at $\tau^{enter}
(t)$, in which case the formula (\ref{pr4})
is not meaningful.
Given a set of  prioritizing functions $\B_{k,p}$, the distribution functions
$U_{kp, ij}^+$ can be determined by setting
\bel{r1}\beta ~\doteq~U^+_{ij}(t) - U^{transit}_{ij}(\tau^{enter}(t))\,,\eeq
\bel{r2}
U_{kp, ij}^+~=~\B_{k,p}(\beta)\,.\eeq
Indeed, $U^+_{ij}(t)$ is the total number of drivers that
reach the end of the arc $\gamma_{ij}$ before time $t$, while
$\beta$ counts how many of these drivers start their journey from the node
$A_i$.
In turn, $\B_{k,p}(\beta)$ determines how many
belong to the subgroup $\G_{k,p}$.

\section{Globally optimal solutions}

In this section we establish the existence of a globally optimal solution.
The proof follows the direct method of the Calculus of Variations,
constructing a minimizing sequence of solutions and showing that a subsequence
converges to the optimal one.
\v
\begin{theorem}\label{theorem1} {\bf (existence of a globally optimal solution).}
Let the flux functions $F_{ij}$ and the cost functions
$\vp_k,\psi_k$ satisfy the assumptions (A1)-(A2).  Then, for
 any $n$-tuple $(G_1,\ldots, G_n)$ of nonnegative numbers,
there exists an admissible  set of departure distributions
$U_{k,p}$ and prioritizing functions $\B_{k,p}$ which yield a
globally optimal solution
of the traffic flow problem.
\end{theorem}
\begin{proof} {\bf 1.}
 By (A2), all functions $\vp_k+\psi_k$ are bounded below.
By possibly adding a constant, it is not restrictive to assume that
$\vp_k(t)+\psi_k(t)\geq 0$ for every time $t$.
Calling $m_0$ the infimum of all total costs in (\ref{JJdef}), taken among
all admissible departure distributions $\{U_{k,p}\}$, this implies
$m_0\geq 0$.   In addition, it is clearly not restrictive to assume
$G_k>0$ for all $k\in \{1,\ldots,n\}$.

Recalling Definitions \ref{definition1} and \ref{definition4}, consider a minimizing
sequence of departure distributions
$U_{k,p}^\nu$ and prioritizing functions $\B_{k,p}^\nu:[0, G_{ij}^\nu]
\mapsto [0, G_{kp}^\nu]$.
Here
$$G_{k,p}^\nu~=~U_{k,p}^\nu(+\infty).$$
By choosing a subsequence, we can assume
\bel{limgk}
\lim_{\nu\to\infty} G^\nu_{k,p}~=~G_{k,p}\quad\qquad\hbox{with}
\qquad\quad \sum_p G_{k,p}~=~G_k\,.\eeq
Moreover,
by Helly's compactness theorem
we can assume that, as $\nu\to\infty$, one has the
pointwise convergence
\bel{pconv}
U_{k,p}^\nu(t)~\to ~U_{k,p}(t)\qquad\qquad \hbox{for a.e.~} t\in\R\,,\eeq
Notice that by \eqref{pf1}, the prioritizing functions $\mathcal{B}^{\nu}_{k,p}$ are Lipschitz continuous with Lipschitz constant 1.  Ascoli's theorem yields the uniform convergence
\bel{uconv}
\B_{k,p}^\nu(\beta)~\to~\B_{k,p}(\beta)\qquad\qquad  \beta\in\R\,.\eeq
In (\ref{uconv}), we have extended the prioritizing
 functions to the entire real line
by setting
$$\B_{k,p}^\nu(\beta)~=~\begin{cases}0\quad &\hbox{if}\quad \beta\leq 0\,,\cr
G_{k,p}^\nu\quad &\hbox{if}\quad \beta\geq G_{ij}^\nu\,.\cr\end{cases}$$
In the remainder of the proof we show that the set of departure distributions
$U_{k,p}$ together with the prioritizing functions $\B_{k,p}$ yield a globally
optimal solution.
\v
{\bf 2.} Let $\ve>0$ be given.   We claim that
there exists a large enough constant $T$, independent of
$\nu$ such that
\bel{tight}\left\{
\bega{rl}\ds\sum_p U_{k,p}^\nu(t)&\leq~\ve\qquad\qquad \hbox{for} ~~t\leq -T\,,\cr
&\cr
\ds\sum_p U_{k,p}^\nu(t)&\geq~G_k-\ve\qquad\qquad \hbox{for} ~~t\geq T\,,\enda
\right.\eeq
for all $\nu$ sufficiently large.
Indeed, by (A2) there exists $T$ such that
$$\vp_k(t)+\psi_k(t)~>~{m_0+1\over\ve}\qquad\qquad
\hbox{for}~~~|t|\geq T\,,~k\in \{1,\ldots,n\}\,.$$
If any one of the two conditions in (\ref{tight}) fails, then the
total cost would be $>m_0+1$.  Since by assumption as $\nu\to\infty$
the total cost approaches the infimum  $m_0$, this proves our claim.

By (\ref{tight}) it follows that the limit functions satisfy
$$U_{k,p}(-\infty) ~=~0,\qquad\qquad U_{k,p}(+\infty) ~=~G_{k,p}\,.$$
In particular, the limit departure distribution is admissible.
\v
{\bf 3.}
For $\beta\in [0, G_{k,p}^\nu]$ let
\bel{betamap}\beta~\mapsto~\tau^{d,\nu}_{k,p}(\beta)\qquad\hbox{and}\qquad
\beta~\mapsto~\tau^{a,\nu}_{k,p}(\beta)\eeq
describe
the departure and arrival time of the $\beta$-driver, in the subgroup $\G_{k,p}$.
We claim that, by possibly extracting a further subsequence, one has
the pointwise convergence
\bel{pwc}\tau^{d,\nu}_{k,p}(\beta)~\to~
\tau^d_{k,p}(\beta)\,,\qquad\qquad \tau^{a,\nu}_{k,p}(\beta)
~\to~\tau^a_{k,p}(\beta)\qquad\qquad\hbox{for a.e.}~\beta\in [0, G_{k,p}].\eeq
Indeed, the maps in (\ref{betamap}) are nondecreasing.  Moreover,
for any $\ve>0$, by (\ref{tight}) these maps
are uniformly bounded when restricted to the subinterval $[\ve , \, G_{k,p}-\ve]$,
uniformly w.r.t.~$\nu$.  By Helly's compactness theorem, we can find a subsequence that
converges pointwise on $[\ve , \, G_{k,p}-\ve]$.   Since $\ve>0$ was arbitrary,
a standard argument proves our claim.
\v
{\bf 4.} It remains to  prove that the limit departure distribution is optimal.
Recall that, without loss of generality, we are assuming
$\vp_k(t)+\psi_k(t)\geq 0$ for all $k,t$.
The total cost $\bar J$ determined by the set of  departure distributions
$\{ U_{k,p}\} $ and prioritizing functions $\{ \B_{k,p}\}$ can be estimated
by
\bel{Jee}\bega{rl}
\bar J
&\doteq\ds ~\sum_{k,p} \int_0^{G_{k,p}} \Big(\vp_k(\tau^d_{k,p}(\beta)) + \psi_k
(\tau^a_{k,p}(\beta))\Big) d\beta\cr
&\cr
&=~\ds \sup_{\ve>0}~\sum_{k,p}\int_\ve^{G_{k,p}-\ve}
\Big(\vp_k(\tau^d_{k,p}(\beta)) + \psi_k
(\tau^a_{k,p}(\beta))\Big) d\beta \,.\enda\eeq
Fix $\ve>0$.  By (\ref{tight}) all departures and arrivals of $\beta$-drivers
with $\beta\in [\ve, \, G_{k,p}-\ve]$ take place in a uniformly bounded
interval of time, say $[-T, T']$.    On this interval, all functions $\vp_k, \psi_k$
are uniformly continuous. Hence the pointwise convergence (\ref{pwc}) yields
$$\bega{rl}&\ds\sum_{k,p} \int_\ve^{G_{k,p}-\ve}
\Big(\vp_k(\tau^d_{k,p}(\beta)) + \psi_k
(\tau^a_{k,p}(\beta))\Big) d\beta\cr &\cr
&\qquad\qquad \ds =~\sum_{k,p} \lim_{\nu\to\infty}\int_\ve^{G_{k,p}-\ve}
\Big(\vp_k(\tau^{d,\nu}_{k,p}(\beta)) + \psi_k
(\tau^{a,\nu}_{k,p}(\beta))\Big) d\beta \cr &\cr
&\qquad\qquad \ds\leq ~\lim_{\nu\to\infty}\sum_{k,p} \int_0^{G_{k,p}}
\Big(\vp_k(\tau^{d,\nu}_{k,p}(\beta)) + \psi_k
(\tau^{a,\nu}_{k,p}(\beta))\Big) d\beta~=~m_0\,.\enda$$
Together with (\ref{Jee}), this  implies $\bar J\leq m_0$,
completing the proof.
\end{proof}

\begin{remark}
{\rm  The above theorem remains valid if in assumption
(A2) we only require   $\vp_k'\leq 0$, $\psi_k'\geq 0$.
}
\end{remark}

\begin{remark} {\rm
A natural conjecture is that, in a globally optimal solution, all departure
rates $\bar u_{k,p}= {d\over dt}U_{k,p}$ are uniformly bounded and have compact support.
 Hence the corresponding
solution should be uniquely determined without need of prioritizing
functions.   In the case of one group of drivers traveling on a single road,
this fact was proved in \cite{BH1}.
In the general case, a proof of the above
conjecture will likely
require a more detailed study of the globally
optimal solution, establishing necessary
conditions for optimality.} \end{remark}
\v

\section{Nash equilibria}

In this section we prove the existence of a Nash equilibrium solution
for traffic flow on a network.  For our model, it turns out that
in a Nash equilibrium all departure rates must be uniformly bounded and have
compact support. As a consequence,  the corresponding
solution is uniquely determined without need of prioritizing
functions.

\v
\begin{theorem}\label{theorem2} {\bf (existence of a Nash equilibrium).}
Let the flux functions $F_{ij}$ and the cost functions
$\vp_k,\psi_k$ satisfy the assumptions (A1)-(A2).
\begi
\item[(i)] For any $n$-tuple
$(G_1,\ldots, G_n)$ of nonnegative numbers, there exists at least one
admissible family of departure rates $\{u_{k,p}^*\}$ which yields a
Nash equilibrium solution.

\item[(ii)] In every Nash equilibrium solution,
all departure rates are uniformly bounded and have
compact support.
\endi
 \end{theorem}
 
\v
Before proving the theorem, we establish a ``modulus of continuity"
for the exit time.   Namely,  for drivers
traveling along a given path $\Gamma_p$,
the arrival time $\tau_p(t)$ is a uniformly continuous function
of the departure time $t$.   In the following, we first consider a single
arc $\gamma_{ij}$  with flux function $F_{ij}(\cdot)$ satisfying (A1).   As in
Figure \ref{f:t94},  $g_{ij}$  denotes the inverse function while
$g^*_{ij}$ is the Legendre transform.
For a driver who enters the arc $\gamma_{ij}$ at time $t$ (possibly joining a queue),
we denote by $\tau_{ij}(t)$ his exit time.

\begin{lemma}\label{Holderlemma} Given  $M>0$, there exists a
continuous function $\phi_{ij}:\R_+\mapsto \R_+$ depending
on the flux function $F_{ij}$ in (\ref{fprop}) and on $M$, such that
$\phi_{ij}(0)=0$ and moreover the following holds.
  Let $U_{ij}^-(\cdot)$ be a Lipschitz
continuous departure distribution,
such that $0\leq u_{ij} (t)={d\over dt} U_{ij}^-(t)\leq M$
for a.e.~$t$.
Then the exit time $\tau_{ij}(\cdot)$ satisfies
\bel{holder}
\tau_{ij}(t_2)-\tau_{ij}(t_1)~\leq~\phi_{ij}(t_2-t_1) \qquad\qquad \hbox{whenever}
\quad  t_1\leq t_2\,.
\eeq
\end{lemma}

\begin{figure}[htbp]
   \centering
 \includegraphics[width=0.9\textwidth]{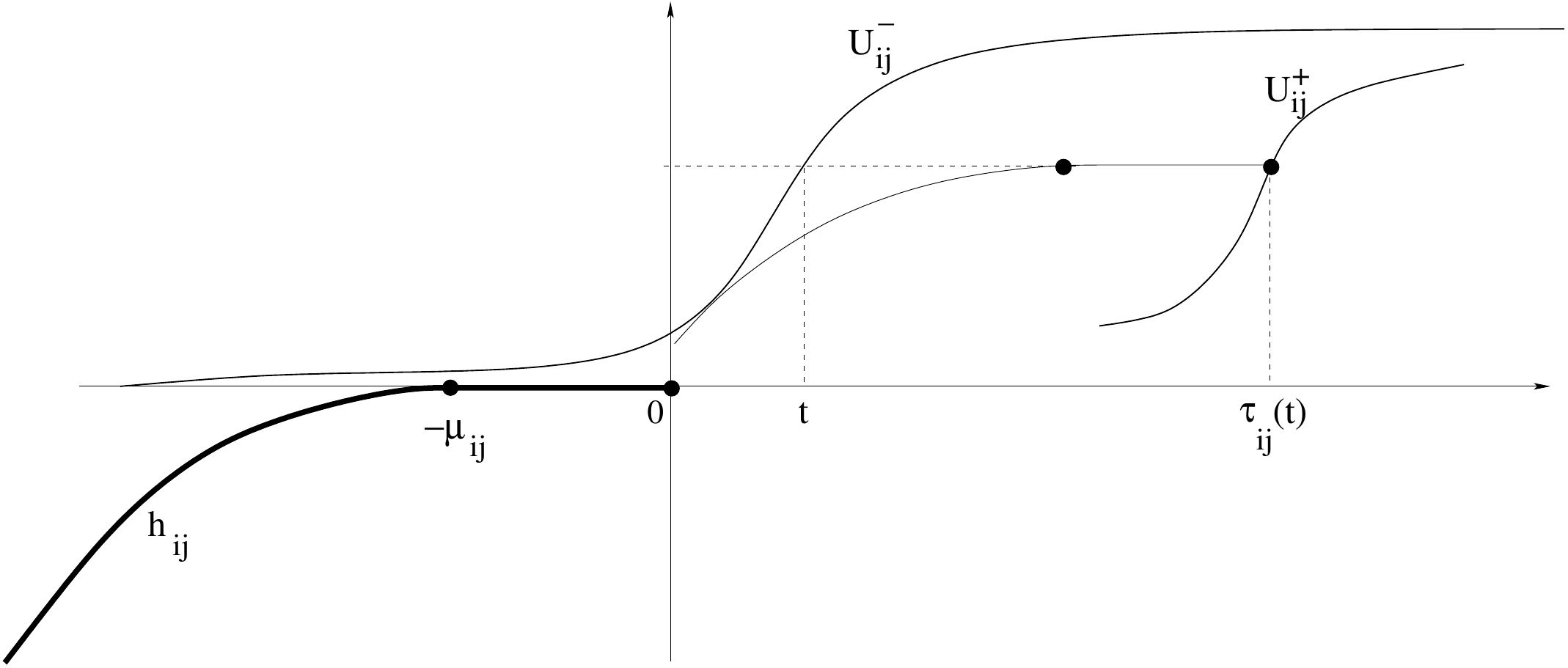}
   \caption{\small  A geometric construction of
   the exit distribution $t\mapsto U_{ij}^+(t)$
   from the entrance distribution $U^-_{ij}$,
   using the formula (\ref{U2}).  }
   \label{f:t105}
\end{figure}

\begin{figure}[htbp]
   \centering
 \includegraphics[width=1.0\textwidth]{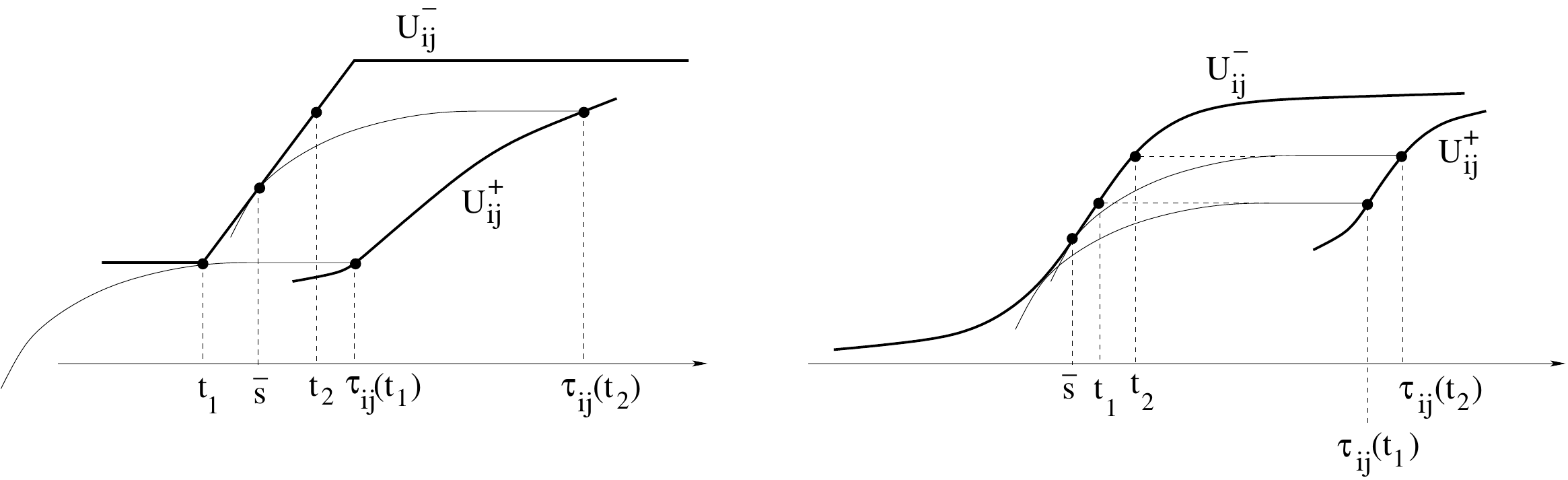}
    \caption{\small  The two cases considered in the proof of Lemma~\ref{Holderlemma}. }
   \label{f:t106}
\end{figure}

\noindent\emph{Proof.} {\bf 1.}   Let $L_{ij}$ be the length of the arc $\gamma_{ij}$.
Following \cite{BH1}, consider the function
\bel{hdef}
h_{ij}(s)~\doteq~-L_{ij}\, g_{ij}^*\left({-s\over L_{ij}}\right)\,.\eeq
Since the Legendre transform $g_{ij}^*$ is convex, the function $h_{ij}$ is concave.   Moreover, calling
$\mu_{ij}\doteq L_{ij}/ v_{ij}(0)$ the minimum time needed to
drive across the arc $\gamma_{ij}$ (= length of the arc divided by the maximum speed),
one has
\bel{hprop} \bega{rl}
h_{ij}(s)&=~ 0   \quad\hbox{if} ~~s\geq -\mu_{ij} \,,\cr
h_{ij}(s)&<~ 0   \quad\hbox{if} ~~s< -\mu_{ij} \,.\enda\eeq
Using the solution formula (\ref{Uij++}), the amount of drivers that exit from the arc $\gamma_{ij}$
before time $\tau$ is computed by (see Figure \ref{f:t105})
\bel{U2}U_{ij}^+(\tau) ~=~\min_s~\Big\{ U_{ij}^-(s)-h_{ij}(s-\tau)\Big\}\,.\eeq
In other words, $U_{ij}^+(\tau)$ is the amount by
which one can shift upward the graph of $h(\cdot - \tau)$,
before hitting the graph of $U^-_{ij}(\cdot)$.
In turn, the exit time of a driver who enters the arc $\gamma_{ij}$ at time $t$
is given by
\bel{xachar}\tau_{ij}(t) ~=~\left( t+ \mu_{ij}\right)~ \vee~
\inf~\Big\{ \tau\,;~~U^-_{ij}(s) \geq U_{ij}^-(t) +
h_{ij}(s-\tau)\quad
\forall s< t\Big\}\,,\eeq
where we used the notation $a\vee b\doteq \max\{a,b\}$.
\v
{\bf 2.}
Let $t_1<t_2$ be given. Two cases can be considered.
\v
\noindent\textbf{Case 1 (Figure \ref{f:t106}, left).} The minimum
\bel{U5}U_{ij}^+(\tau_{ij}(t_2)) ~=~\min_s~\Big\{ U_{ij}^-(s)-h_{ij}(s-\tau_{ij}(t_2)
)\Big\}\,.\eeq
is attained at a point $\bar s\in [t_1, t_2]$.
In this case we have
$$\tau_{ij}(t_1)~\geq~t_1+\mu_{ij}\,,$$
$$\tau_{ij}(t_2)~\leq~t_2 + \inf\Big\{ \tau\,;~~Ms\geq
h_{ij}(s-\tau)\quad
\forall s\in [- M(t_2-t_1)\,,~ 0]\Big\}\,.$$
Hence
\bel{bb1}
\tau_{ij}(t_2)-\tau_{ij}(t_1)~\leq~\phi^\sharp_{ij}(t_2-t_1)\eeq
where the continuous increasing function $\phi^\sharp_{ij}$ is defined by
\bel{bb2}
\phi^\sharp_{ij}(\xi)~\doteq~\xi+
\inf\Big\{ \tau\,;~~Ms\geq
h_{ij}(s-\tau)\quad
\forall s\in [- M\xi\,,~ 0]\Big\} - \mu_{ij}\,.\eeq
Notice that (\ref{bb1})-(\ref{bb2}) estimate
a ``worst case scenario" where the driver
starting at time $t_1$ finds an empty road ahead, while the driver starting
at time $t_2>t_1$ has $M(t_2-t_1)$ cars in front of him.
\v
\noindent\textbf{Case 2 (Figure \ref{f:t106}, right).}
The minimum in (\ref{U5}) is attained at a point $\bar s <t_1$.
In this case we have
$$U_{ij}^-(t_1) +h_{ij}(\bar s- \tau_{ij}(t_1)) ~\leq~U_{ij}^-(\bar s)~
~=~U^-_{ij}(t_2)+h_{ij}(\bar s-\tau_{ij}(t_2))\,.$$
Since $h_{ij}$ is concave and $\tau_{ij}(t_1)- t_1\geq \mu_{ij}$,
~ $\tau_{ij}(t_2)- t_2\geq \mu_{ij}$, one has
$$\bega{rl}0&\geq~h_{ij}\left(\tau_{ij}(t_1)-\tau_{ij}(t_2)-\mu_{ij}\right)~\geq~
h_{ij}(\bar s- \tau_{ij}(t_2))-h_{ij}(\bar s- \tau_{ij}(t_1)) \cr
&\cr
&\geq~U_{ij}^-(t_1)-U_{ij}^- (t_2)~\geq~ M (t_1-t_2).\enda
$$
In this case, we conclude
\bel{bb3}
\tau_{ij}(t_2)-\tau_{ij}(t_1)~\leq~\phi^\flat_{ij}(t_2-t_1),\eeq
where the continuous function $\phi_{ij}^\flat$ is implicitly defined by
\bel{bb4}
h_{ij}(-\mu_{ij}-\phi^\flat_{ij}(\xi))~=~-M\xi\,.\eeq

{\bf 3.} The two functions $\phi_{ij}^\sharp$, $\phi_{ij}^\flat$
defined at (\ref{bb2}), (\ref{bb4}) are both continuous
and vanish at the origin.
Defining
\bel{bb5}\phi_{ij}(\xi)~\doteq~\max~\{\phi^\sharp_{ij}(\xi)\,,~\phi^\flat_{ij}(\xi)
\}\,,\eeq
the conclusion of the lemma is satisfied.
\endproof
\v
The next lemma extends the result of Lemma \ref{Holderlemma} to a general
path $\Gamma_p$. Here  $\tau_p(t)$ denotes
the arrival time of a driver
starting at time $t$ and traveling along $\Gamma_p$.

\begin{lemma}\label{networklemma} Let all departure rates  $\ov u_{k,p}(t)$
be uniformly bounded, so that
\bel{unifbdd}
\ov u_{k,p}(t)~\leq~M_0 \qquad \forall ~~k,p, t\,.
\eeq

Then, for
any viable path $\Gamma_p$, there exists continuous function $\phi_p:\R_+\mapsto \R_+$
such that
\bel{holdertt}\phi_p(0)~=~0\,,\qquad\qquad
\big(\tau_{p}(t_2)-\tau_{p}(t_1)\big)~\leq~\phi_p(t_2-t_1)
\qquad\qquad \hbox{whenever}
\quad t_1\leq t_2 \,.
\eeq
\end{lemma}

\noindent\emph{Proof.} The assumption (\ref{unifbdd}), together with the fact that the flux
through
each arc $\gamma_{ij}$ cannot be greater than the
maximum flux $F_{ij}^{max}$, implies
that the maximum incoming flux through each arc is
bounded by some constant $M$.
For each arc $\gamma_{ij}$, let
$\phi_{ij}$ be the modulus of continuous
dependence constructed in Lemma~\ref{Holderlemma}.
If $\Gamma$ is the path  in (\ref{path}), obtained as the concatenation of the arcs
$\gamma_{i(\ell-1), i(\ell)}$, $\ell = 1,\ldots,\nu$, it now suffices to define the
function $\phi_p$ as the composition of the corresponding functions
$\phi_{i(\ell-1), i(\ell)}$, namely
$$\phi_p~\doteq~\phi_{i(\nu-1), i(\nu)}\circ\cdots\circ \phi_{i(1), i(2)}
\circ
\phi_{i(0), i(1)}\,.$$
\endproof

Our final lemma shows that the arrival times, and hence the cost functions,
depend continuously on the
departure rates.
\begin{lemma}\label{clem} Consider a sequence of departure rates  $u^\nu =
(u_{k,p}^\nu)$ which are
uniformly bounded and  supported inside a common interval $I= [-T,T]$, namely
\bel{ubbb}
0~\leq~ u^\nu_{k,p}(t)~\leq~M_0 \qquad \forall ~~t\in I\,,\qquad\qquad
u^\nu_{k,p}(t) ~=~0\qquad \hbox{if}~~ t\notin I\,,
\eeq
for every $\nu\geq 1$.
Assume that for all $k,p$ one has  the weak convergence in $\L^1$
\bel{wco}u^\nu_{k,p}~\wto~  u^*_{k,p}\,.\eeq
For each viable path $\Gamma_q$, call $ \tau_q^\nu(t)$,
$\tau^*_q(t)$  the corresponding arrival times
of a driver who departs at time $t$ and travels along $\Gamma_q$.
Then, as $\nu\to\infty$, one has the uniform convergence
\bel{tconv}
\|\tau_q^\nu - \tau^*_q\|_{\C([-T,T])}~\to ~0\,.\eeq
\end{lemma}
\v
\noindent\emph{Proof.} {\bf 1.}
We first consider the case of a single arc $\gamma=\gamma_{ij}$.
Assume that the departure rates $u^\nu$ satisfy
$$\bega{cl}u^\nu(t)~\in ~[0,M_0]\qquad &\hbox{if}~~~t\in[-T,T]\,,
\cr u^\nu(t)~=~0\qquad &\hbox{if}~~~t\notin[-T,T]\,,\enda$$ and converge weakly:
$u^\nu\rightharpoonup  u$.
Define the integrated functions
\begin{equation}\label{uandU}
U^\nu(t)~\doteq~\int_{-\infty}^tu^\nu(s)\,ds,\qquad\qquad
 U(t)~\doteq~\int_{-\infty}^t  u(s)\,ds\,.
\end{equation}
Our assumptions imply the convergence
$U^\nu(t)\rightarrow  U(t)$, uniformly  for $t\in \R$.
In turn,  this implies that the exit distributions
\bel{exitn}
U^{\nu}_+(t)\doteq \min_{\tau}\left\{U^{\nu}(\tau)
+L_{ij}\,g_{ij}^*\Big({t-\tau\over L_{ij}}\Big)\right\},\,
  U_+(t)\doteq\min_{\tau}\left\{ U(\tau)+L_{ij}
  \,g_{ij}^*\Big({t-\tau\over L_{ij}}\Big)\right\},\,
\eeq
satisfy the uniform convergence $U^\nu_+\to U_+$.
Indeed, by (\ref{exitn}) we have
\bel{uncon}
\|U^\nu_+- U_+\|_{\L^\infty}~\leq~
\|U^\nu- U\|_{\L^\infty}~\to ~0\qquad\qquad \hbox{as}~~\nu\to\infty\,.\eeq

\v
{\bf 2.} Given $\ve>0$ choose $\hat \nu $ large enough so that
\bel{Unu}|U^\nu(t)- U(t)|~\leq~\ve\qquad\qquad \forall \nu\geq \hat \nu \,,~~t\in\R\,.
\eeq
Fix a time $t\in [-T,T]$. To estimate the difference $|\tau^\nu(t)-\tau(t)|$
between the corresponding
arrival times, consider the modified departure distribution
\bel{Ue}
U^\sharp(s)~\doteq~\left\{ \bega{cl}  U(s)\quad &\hbox{if}~~s\leq t\,,\cr
 U(t) + (2+M_0)(s-t)\quad &\hbox{if}~~s\in [t,T]\,,\cr
 U(t) + (2+M_0)(T-t)\quad &\hbox{if}~~s\geq T\,.\cr
\enda\right.\eeq
Call $\tau^\sharp(t)$ the arrival time of a driver departing at time $t$,
relative to
the distribution $U^\sharp$.
By Lemma~\ref{Holderlemma}, the function $t\mapsto\tau^\sharp(t)$ satisfies
a uniform modulus of continuity say
$$\tau^\sharp(t_2)-\tau^\sharp(t_1)~\leq~\phi(t_2-t_1)\qquad\qquad\forall t_1<t_2\,,
$$
for some continuous function $\phi$ with $\phi(0)=0$.
In particular,
\bel{tee}
\tau^\sharp(t+\ve)-\tau^\sharp(t) ~\leq~\phi(\ve)\,.\eeq
Observing that
$$U^\sharp(t+\ve)-U^\sharp(s)~\geq~U^\nu(t+\ve)-U^\nu(s)\qquad\qquad \forall
\nu\geq \hat \nu \,,~~s\leq t+\ve$$
and using the representation formula (\ref{xachar})
with $h_{ij}$ given by (\ref{hdef}),
we obtain
$$\bega{rl}&\tau^\sharp(t+\ve)\cr &\cr
 =&(t+\ve+\mu_{ij})\vee
\inf\left\{ \tau\,;~U^\sharp(t+\ve) -U^\sharp(s)
\leq L_{ij}\, g_{ij}^*\Big({-s+\tau\over L_{ij}}\Big)~~\hbox{for all}~s\leq
t+\ve\right\}\cr
&\cr
\geq &(t+\ve+\mu_{ij})\vee \inf\left\{ \tau\,;~U^\nu(t+\ve) -U^\nu(s)
\leq L_{ij}\, g_{ij}^*\Big({-s+\tau\over L_{ij}}\Big)~~\hbox{for all}~s\leq
t+\ve\right\}\cr
&\cr
 =&\tau^\nu(t+\ve)\,.\enda$$
Therefore
\bel{in5}\tau^\nu(t)~\leq~\tau^\nu(t+\ve)~\leq~\tau^\sharp(t+\ve)~\leq~\tau^\sharp
(t)+\phi(\ve)
~=~\tau(t)+\phi(\ve) \qquad\forall
\nu\geq \hat \nu \,.\eeq
Switching the roles of $U,U^\nu$ and hence defining
$$
U^\sharp(s)~\doteq~\left\{ \bega{cl} U^\nu(s)\quad &\hbox{if}~~s\leq t\,,\cr
U^\nu(t) + (2+M_0)(s-t)\quad &\hbox{if}~~s\in [t,T]\,,\cr
 U^\nu(t) + (2+M_0)(T-t)\quad &\hbox{if}~~s\geq T\,,\cr
\enda\right.$$
we obtain the converse inequality
\bel{in6}\tau(t)~\leq~\tau^\nu(t)+\phi(\ve)\qquad\qquad\forall
\nu\geq \hat \nu \,.\eeq
\v
{\bf 3.} For each arc $\gamma_{ij}$,  the inequalities (\ref{in5})-(\ref{in6}) show that the functions $t\mapsto \tau_{ij}^\nu(t)$
converge uniformly to the corresponding functions $t\mapsto \tau_{ij}(t)$.
Given any path $\Gamma_q$, the uniform convergence $\tau_q^\nu\to \tau^*_q$
is now obtained by a straightforward argument, using
induction on the number of arcs contained in $\Gamma_q$.
\endproof

\v
We are now ready to prove the main result of this paper.

\noindent\emph{Proof of Theorem \ref{theorem2}.}
{\bf 1.}
We claim that there exists a time interval $I=[-T_0,\,T_0]$  so large that, in any Nash equilibrium, no driver will depart or arrive at a time $t\notin I$.
Indeed, given the $n$-tuple $(G_1,\ldots, G_n)$, the travel time along
any viable path $\Gamma_p=\big(\gamma_{i(0),i(1)},\ldots,
\gamma_{i(\nu-1),i(\nu)}\big)$ is a priori bounded by
\bel{maxtt}
T^{max}_p~\doteq~\sum_{\ell=1}^{\nu}\left\{
{G\over F_{i(\ell-1),i(\ell)}^{max}}+
{L_{i(\ell-1),i(\ell)}\over v_{i(\ell-1),i(\ell)}
\big(\rho^*_{i(\ell-1),i(\ell)}\big)}\right\}.
\eeq
Here and in the sequel, we call
\bel{Gdef}
G~\doteq~G_1+\cdots+G_n\eeq
the total number of drivers.
Notice that, in each summand on the right hand side of (\ref{maxtt}), the first
term is an upper bound for the time spent waiting in the queue (total number of drivers divided by the maximum flux)
while the second term
is an upper bound on the actual travel time (length divided by the minimum speed).
Let
$$T^{max}~\doteq ~\max_p ~T^{max}_p$$
be an upper bound on the travel time along all viable paths.
 In view of assumption  (A2), there exists $T_0$ large enough such that
\bel{Idef}
\min_{k}\big\{\vp_k(t)+\psi_k(t)\big\}~>~
\max_k\big\{\vp_k(0)+\psi_k(T^{max})\big\}\quad \forall t\notin I \doteq [-T_0,\,T_0]\,.
\eeq
Therefore, in a Nash equilibrium no driver will depart or arrive outside $[-T_0,\,T_0]$.
Otherwise, he would achieve a strictly lower cost by departing at time $t=0$.
\v
{\bf 2.}  Let $F^{max}=\ds \max_{i,j} F^{max}_{ij}$
be an upper bound for the flux over all arcs.
Call
$$\vp'_{max}~\doteq~\max_{1\leq k\leq n}
\max_{t\in I} ~|\vp_k'(t)|\,,\qquad\qquad \psi'_{min}
~\doteq~\min_{1\leq k\leq n} \min_{t\in I} ~\psi_k'(t).$$
Observe that $\psi'_{min}>0$, because of the assumption (A2).
We claim that, in a Nash equilibrium, all departure rates $u_{k,p}$ must satisfy the
 priori bound
\bel{kdef}
u_{k,p}(t)~\leq~\kappa~\doteq~{\vp'_{max}\cdot F_{max}\over \psi'_{min}}
\qquad\qquad\hbox{for a.e.~}t\,.\eeq
Indeed, consider drivers of the $k$-th family traveling along the path $\Gamma_p$.
Let $t_1<t_2$ be any two departure times, and call $\tau_1<\tau_2$
the corresponding arrival times.
The total costs for these two drivers must be the same, hence
$$\vp(t_1)+\psi(\tau_1) ~=~\vp(t_2)+\psi(\tau_2)\,.$$
On the other hand, the upper bound on the flux implies
$$\tau_2-\tau_1 ~\geq~{1\over F_{max}}\int_{t_1}^{t_2} u_{k,p}(t)\, dt\,.$$
Therefore
\begin{multline}
(t_2-t_1)\,\vp'_{max}
~\geq~\vp(t_1)-\vp(t_2) ~=~\psi(\tau_2)-\psi(\tau_1)
\\
~\geq~(\tau_2-\tau_1)\,\psi'_{min}
~\geq~{\psi'_{min}\over F_{max}}\int_{t_1}^{t_2} u_{k,p}(t)\, dt\,.
\end{multline}
We thus conclude
$${\vp'_{max}\cdot F_{max}\over \psi'_{min}} ~\geq ~{1\over  t_2-t_1} \int_{t_1}^{t_2}
u_{k,p}(t)\, dt\,.$$
Since this bound is valid for every interval $[t_1, t_2]\subseteq I$,
the pointwise bound (\ref{kdef}) must hold.  Moreover,
for $t\notin I$ we already know that $u_{k,p}(t)=0$.
The last statement of the Theorem is thus proved.
 \v
{\bf 3.}  Choose the time
\bel{Tdf}
T~\doteq~T_0 + {G\over \kappa}\,,\eeq
where $G\doteq \sum_{k}G_k \,$.

Consider the family of  admissible departure rates
\bel{Udep}\bega{rl}
\U&\doteq~\bigg\{ (u_{k,p})_{1\leq k\leq n,\, 1\leq p\leq N}~;\quad
u_{k,p}:\R\mapsto [0,4\kappa]\,,\quad u_{k,p}(t) =0 ~~\hbox{for}~~t\notin [-T,T]\,,
\cr
&\cr
&\ds\qquad \qquad
\sum_p \int  u_{k,p}(t)\, dt ~=~G_k\quad\hbox{for every}~k\bigg\}.\enda\eeq
It is understood that $u_{k,p}\equiv 0$ if the path $\Gamma_p$ does not connect
$A_{d(k)}$ with $A_{a(k)}$.
Notice that $\U$ is a closed convex subset of $\L^1(\R;\,\R^{n\times N})$.

For each fixed $\nu\geq 1$, we consider a finite dimensional subset $\U_\nu\subset\U$
consisting of all $u=(u_{k,p})$ which are piecewise constant on time
intervals of length
$T/\nu$.
Introducing the points
$$t_\ell^\nu ~\doteq~ {\ell \over \nu} \,T \,,\qquad\qquad -\nu\leq \ell\leq \nu\,,$$
we thus define
\bel{unudef}\bega{rl}\U_\nu&\doteq~\ds
\Big\{ u=(u_{k,p})\in \U\,;
\cr &\qquad \ds \hbox{every function  $u_{k,p}$
is constant on each subinterval
$I^\nu_\ell~\doteq~\,]t_{\ell-1}^\nu,\, t_\ell^\nu]$}\Big\}.\enda\eeq
Observe that  every  $u\in \U_\nu$ has the form
\bel{ukpl}
u = (u_{k,p})\,,\qquad\qquad u_{k,p}(t) ~=~u_{k,p,\ell}\qquad\forall \quad t\in \,]t_{\ell-1}^\nu,\, t_\ell^\nu]\,,
\eeq
for some constants $u_{k,p,\ell}\in [0, \, 4\kappa]$.
\v
{\bf 4.} Given $u=(u_{k,p})\in \U$, let $\tau_{q}(t)$ be the arrival time
of a driver starting at time $t$ and traveling along the path $\Gamma_q$.
Clearly, this arrival time depends on the overall traffic conditions, hence on
all functions $u_{k,p}$. If this driver belongs to the $j$-th family, his total
cost is
$$\Phi^{(u)}_{j,q}(t) ~=~\vp_j(t) + \psi_j(\tau_{q}(t))\,.$$
We now observe that, for each $\nu\geq 1$, the domain
$\U_\nu$ is a finite dimensional, compact, convex
subset of $\L^2([-T,T];\, \R^{n\times N})$.  Moreover, by Lemma~\ref{clem} the maps
$u\mapsto \Phi^{(u)}_{k,p}(\cdot)$ are continuous from $\U_\nu$ into $\L^2$.
Hence, by the  theory of variational inequalities \cite{KS},
there exists a function  $\bar u^\nu=(\bar u^\nu_{j,q})\in \U_\nu$ which
satisfies
\bel{vi1}
\sum_{j,q}\int_{-T}^T \Phi_{j,q}^{(\bar u^\nu)}(t) \cdot \Big( v_{j,q}(t)-
\bar u_{j,q}^\nu(t)\Big)\, dt~\geq~0\qquad\qquad \forall v\in \U_\nu\,.\eeq
\v
{\bf 5.} We now let $\nu\to\infty$.
By the previous steps, there exists a sequence of
piecewise constant
functions $\bar u^\nu= (\bar u^\nu_{k,p}) \in \U_\nu$ such that (\ref{vi1}) holds
for every $\nu\geq 1$.
Since all functions $\bar u^\nu_{k,p}$ are uniformly bounded and supported inside the
interval $I = [-T,T]$,
by taking a subsequence we can assume the weak convergence
\bel{wlim}(\bar u^\nu_{k,p})~\wto~(u^*_{k,p})\eeq
for some function $u^*= (u^*_{k,p})\in \U$.  We claim that the departure rates
$u^*_{k,p}$
yield a Nash equilibrium solution. More precisely:

\begi\item[(NE)]
Given any $k,p$, any  $t_1\in Supp(u^*_{k,p})$,  $t_2 \in \R$ and any
path $\Gamma_q$ with the same initial and final nodes as $\Gamma_p$,
one has
\bel{samec}\Phi^*_{k,p}(t_1)~\leq ~\Phi^*_{k,q}(t_2).\eeq
\endi

Indeed,  (\ref{samec}) implies that no $k$-driver can lower his own cost
by switching to the time $t_2$ or choosing the alternative path $\Gamma_q$ to reach destination.  We recall that $t$ is in the support of
a function $f\in \L^1$ if and only if  $\int_{t-\ve}^{t+\ve}|f(s)|\, ds \not= 0$
for every $\ve>0$.
\v
{\bf 6.}
By Lemma \ref{networklemma}, all maps $t\mapsto \tau^{\nu}_{k,p}(t)$ have the same modulus
of continuity.  Since these arrival times are uniformly bounded, we can apply
 Ascoli's compactness theorem.
Choosing a subsequence and relabeling, as $\nu\to\infty$   we
thus   achieve the  convergence
\bel{tconv2}
\tau^{\nu}_{k,p}(t)~\to~\tau^*_{k,p}(t)\qquad\qquad \hbox{for all $ k,p$,
uniformly for $t\in [-T,T]$.}
\eeq
By the assumption {\bf (A2}), the cost functions $\vp_k(\cdot),\,\psi_k(\cdot)$ are
 continuous.  Therefore,
 the functions $\Phi^{\nu}_{k,p}(\cdot)=\vp_k(\cdot)+\psi_k\big(\tau_{k,p}^{\nu}(\cdot)\big)$
 converge to  $ \Phi^*_{k,p}(\cdot)=\vp_k(\cdot)+\psi_k\big(\tau_{k,p}^*(\cdot)\big)$,
uniformly for $t\in [-T,T]$.  By Lemma~\ref{clem}, $\tau_{k,p}^*(t)$ is indeed the
arrival time of a $k$-driver departing at time $t$ and following the path $\Gamma_p$,
in the case where the departure rates of all drivers are given by $u^*= (u^*_{j,q})$.
\v

\begin{figure}[htbp]
 \c{\includegraphics[width=1.0\textwidth]{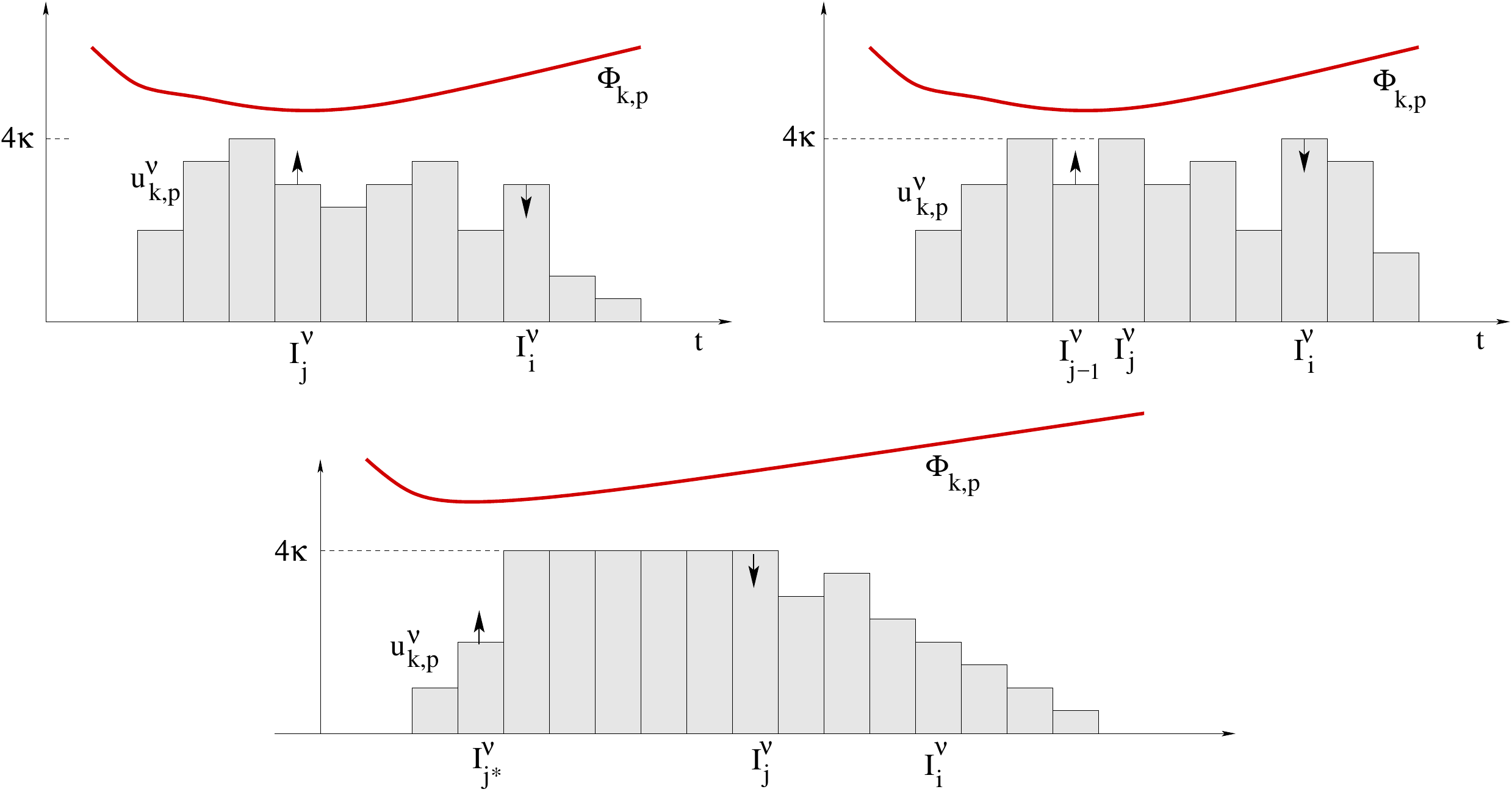}}
    \caption{\small  The three cases considered in the proof of Theorem \ref{theorem2}.
    In Case 1 (top left) the average of the cost
    $\Phi_{k,p}$ on the interval $I_i^\nu$ is higher than
    on the interval $I_j^\nu$.  To obtain a contradiction with
     (\ref{vi1}) we simply move some of the mass from
    $I^\nu_i$ to $I^\nu_j$.   In Case 2a (top right) one cannot increase the value
    of $u^\nu_{k,p}$ on the interval $I_j^\nu$ because of the constraint
    $u\leq 4\kappa$.   However, some mass can be moved from $I_i^\nu$
    to the previous interval
    $I_{j-1}^\nu$.    In Case 2b (bottom) there are several adjacent intervals
    where $u^\nu_{k,p}\equiv 4\kappa$.    In this case, if $I^\nu_{j^*}$ is the first interval to the left of $I^\nu_j$ where $u^\nu_{k,p}< 4\kappa$, we argue that
  (i)   $t_{j^*}^\nu> -T$, and (ii)  the average of the cost
     $\Phi_{k,p}$  on $I^\nu_{j^*}$ is strictly less
     than on $I^\nu_j$.  In this last case, to obtain a contradiction with
     (\ref{vi1}) we move some mass from
     $I^\nu_{j^*}$ to $I^\nu_j$.
      }
   \label{f:tf19}
\end{figure}
{\bf 7.} If
 (\ref{samec}) fails,
then by continuity there exists $\delta>0$ such that
\bel{diffc}\Phi^*_{k,p}(t)~>~\Phi^*_{k,q}(t')+2\delta
\qquad\qquad\hbox{whenever} ~|t-t_1|\leq 2\delta, ~~|t'-t_2|\leq 2\delta\,.\eeq
By uniform convergence, for all $\nu$ large enough we have
\bel{diffc2}\Phi^\nu_{k,p}(t)~>~\Phi^\nu_{k,q}(t')+\delta
\qquad\qquad\hbox{whenever} ~|t-t_1|\leq 2\delta, ~~|t'-t_2|\leq 2\delta\,.\eeq
Observe that it is not restrictive to assume
that $t_2\in [-T_0,T_0]$.
Indeed, if (\ref{samec}) fails for some $t_2\notin [-T_0,T_0]$, then (\ref{Idef})
implies
$$
\Phi^*_{k,p}(t_1)~> ~\Phi^*_{k,q}(t_2) ~>~ \Phi^*_{k,q}(0),$$
and we can simply replace $t_2$ with zero.

The weak convergence (\ref{wlim}), together with the assumption on the support
of the function $u^*_{k,p}$,  now  implies
$$\lim_{\nu\to\infty} \int_{t_1-\delta}^{t_1+\delta} \bar u_{k,p}^\nu(t)\, dt~=~ \int_{t_1-\delta}^{t_1+\delta} u_{k,p}^*(t)\, dt~>~0\,.$$
Therefore, for every $\nu$ sufficiently large we can find two intervals
\bel{in} I^{\nu}_{i}~=~]t^\nu_{i-1},\,t^\nu_{i}]~\subset~ [t_1-\delta\,,~t_1+\delta],\qquad\qquad
I^\nu_j~=~]t^\nu_{j-1},\,t^\nu_j]~\subset~ [t_2-\delta\,,~t_2+\delta]\eeq
with $t_j^\nu>t_2$ and $\bar u^\nu_{k,p}(t)>0$ for  $t\in I^\nu_i$.
\v

{\bf 8.} We now derive a contradiction, showing that, for $\nu$ sufficiently large,
the departure rates $\bar u^\nu_{k,p}$ do not satisfy the variational inequality
(\ref{vi1}).
Two possibilities can arise (see Figure \ref{f:tf19}).
\v
\noindent\textbf{Case 1.} $\bar u^\nu_{k,q,j} ~<~4\kappa$.   In this case we  define
a new set of departure rates $v^\ve = ( v^\ve_{k,p})\in \U_\nu$
by setting
$$v^\ve_{k,p}(t) ~=~u^\nu_{k,p}(t) - \ve\qquad \hbox{if}\quad t\in I^\nu_i\,,$$
$$v^\ve_{k,q}(t) ~=~u^\nu_{k,q}(t) + \ve\qquad\hbox{if}\quad  t\in I^\nu_j\,,$$
and setting $v^\ve_{h,r}(t) = u^\nu_{h,r}(t)$ in all other cases.
Notice that, if $\ve=\min\{ u^\nu_{k,p,i}, ~4\kappa- u^\nu_{k,q,j}\}$
then these new departure rates are still admissible.
By (\ref{diffc2}) and (\ref{in}), this construction yields
\begin{multline}\label{vi2}
\sum_{h,r}\int_{-T}^T \Phi_{h,r}^{(\bar u^\nu)}(t) \cdot \Big( v^\ve_{h,r}(t)-
\bar u_{h,r}^\nu(t)\Big)\, dt
\\
~=~\ve\int_{I^\nu_j} \Phi_{k,q}^{(\bar u^\nu)}(t)\, dt - \ve
\int_{I^\nu_i} \Phi_{k,p}^{(\bar u^\nu)}(t)\, dt~\leq~ -2\ve \delta\,,
\end{multline}
providing a contradiction with (\ref{vi1}).
\v
\noindent\textbf{Case 2.} $\bar u^\nu_{k,q,j} ~=~4\kappa$. If this equality holds, consider the
index
$$j^*~\doteq~\max ~\{ i<j\,;~~\bar u^\nu_{k,q,i} ~<~4\kappa\}\,.$$
Notice that $t^\nu_{j^*}>-T$.
Indeed, by construction $t_2>-T_0$. If $t^\nu_{j^*}\leq -T$, by (\ref{Tdf})
this would imply
$$\bar u^\nu_{k,q}(t) = ~4\kappa\qquad\qquad \forall ~t\in  [t^\nu_{j^*}\,, t^\nu_j]
~\supseteq~[-T, \, -T_0]\,,$$
$$\int \bar u^\nu_{k,q}(t)\, dt~\geq~4\kappa (t^\nu_j - t^\nu_{j^*})~
~\geq~4\kappa(T-T_0)>~G\,,$$
reaching a contradiction.
We consider two  subcases.

\noindent\textbf{Case 2a.} $j^* = j-1$.   In this case, since it is not restrictive to assume
${T\over \nu}<{\delta\over 4}$, we have $I^\nu_{j-1}
=[t^\nu_{j-2}, t^\nu_{j-1}]\subset[t_2-\delta,~t_2+\delta]$.
We can thus derive a contradiction as in Case 1, simply replacing $j$ by $j-1$.
\v
\noindent\textbf{Case 2b.} $j^* \leq j-2$.  Observe that, for all $s_1<s_2$,
\bel{fifo}
\tau_{k,q}^{\nu}(s_2)-\tau_{k,q}^{\nu}(s_1)~\geq~{1\over F_{max}}\int_{s_1}^{s_2}
\bar u^\nu_{k,q}(\xi)\, d\xi\,.
\eeq
In particular, for any
$s_1\in I^\nu_{j^*}$ and $s_2\in I^\nu_j$ we have
 \begin{multline}\label{fi5}
\tau_{k,q}^{\nu}(s_2)-\tau_{k,q}^{\nu}(s_1)~\geq~{1\over F_{max}}\int_{s_1}^{s_2}
\bar u^\nu_{k,q}(\xi)\, d\xi
\\
~\geq~{1\over F_{max}} \, 4\kappa\,[ t_{j-1}-t_{j^*}]
~\geq~{4\kappa (s_2-s_1)\over 3 F_{max}}.
\end{multline}
This yields the estimate
\begin{multline}\label{est1}
\psi_{k}\big(\tau^{\nu}_{k,q}(s_2)\big)-\psi_k\big(\tau^{\nu}_{k,q}(s_1)\big)
\\
~\geq~\psi'_{min}\big(\tau^{\nu}_{k,q}(s_2)-
\tau^{\nu}_{k,q}(s_1)\big)~\geq~\psi'_{min}\cdot {4\kappa (s_2-s_1)\over 3 F_{max}}\,.
\end{multline}
On the other hand, we have
\bel{est2}
\vp_k(s_2)-\vp_k(s_1)~\geq~-\vp'_{max}(s_2-s_1)\,.
\eeq
Recalling the definition of the constant $\kappa$ in (\ref{kdef}),
from (\ref{est1})-(\ref{est2}) we obtain
\bel{ett}
\Phi_{k,q}^{\nu}(s_2)-\Phi_{k,q}^{\nu}(s_1)~\geq~
\left({4\kappa\psi'_{min}\over3 F_{max}}-\vp'_{max}\right)(s_2-s_1)~
=~{1\over 3}\vp'_{max}\cdot (s_2-s_1)\eeq
for all $s_1\in I^\nu_{j^*}$ and $s_2\in I^\nu_j\,$.

We now choose the
departure rates $v^\ve = ( v^\ve_{k,p})\in \U_\nu$
by setting
\bel{vep2}
\bega{c}v^\ve_{k,p}(t) ~=~u^\nu_{k,q}(t) - \ve\qquad \hbox{if}\quad t\in I^\nu_j\,,
\cr v^\ve_{k,q}(t) ~=~u^\nu_{k,q}(t) + \ve\qquad\hbox{if}\quad  t\in I^\nu_{j^*}\,,
\enda\eeq
and setting $v^\ve_{h,r}(t) = u^\nu_{h,r}(t)$ in all other cases.
Notice that, if $\ve=\min\{ u^\nu_{k,q,j}, ~4\kappa- u^\nu_{k,q,j^*}\}$
then these new departure rates are still admissible.

Using (\ref{ett}) with $s_1=t$, $s_2=t+t^\nu_j - t^\nu_{j^*}$ we compute
\bel{vi3}\bega{rl}
&\ds\sum_{h,r}\int_{-T}^T \Phi_{h,r}^{(\bar u^\nu)}(t) \cdot \Big( v^\ve_{h,r}(t)-
\bar u_{h,r}^\nu(t)\Big) dt~
=~\ve\left(\int_{I^\nu_{j^*}} \Phi_{k,q}^{(\bar u^\nu)}(t)dt -  \int_{I^\nu_j}
\Phi_{k,q}^{(\bar u^\nu)}(t)  dt\right)\cr &\cr
&=~\ds \ve
\int_{I^\nu_{j^*}} \Big(\Phi_{k,q}^{(\bar u^\nu)}(t) -
 \Phi_{k,q}^{(\bar u^\nu)}(t+t^\nu_j - t^\nu_{j^*})\Big)\,
  dt~\leq~ -{\ve\over 3}\vp'_{max}\cdot (t^\nu_j - t^\nu_{j^*})
~<~0\,.\enda\eeq
Once again we reached a contradiction with (\ref{vi1}), completing the proof.
\endproof
\v
\begin{remark} {\rm The above proof is based on a topological method, and does
not yield any information about the uniqueness of the Nash equilibrium.
Another important question is the dynamic stability of this equilibrium solution.
This issue was investigated
numerically in \cite{BLSY}. In the case of
a single group of drivers traveling on a single road, the uniqueness
of the Nash equilibrium solution was proved in \cite{BH1},
but the stability issue remains
unresolved.}\end{remark}



\medskip
Received August 2012; revised February 2013.
\medskip

\end{document}